\documentclass[letterpaper, 10pt, conference]{ieeeconf}
\IEEEoverridecommandlockouts        
\usepackage{hyperref}

\title{\LARGE \bf Mode Participation and Inter-Area-Observability Blocking Controllers for Power Networks }

\author{Rajasekhar Anguluri and Abdullah Al Maruf
	\thanks{R.Anguluri is with the Department of Computer Science and Electrical Engineering, University of Maryland, Baltimore County, MD 85281, USA. A. Al Maruf is with the Department of Electrical and Computer Engineering, California State University, Los Angeles, CA 90032, USA (emails: \href{rajangul@umbc.edu}{rajangul@umbc.edu} and \href{amaruf@calstatela.edu}{amaruf@calstatela.edu}). The authors thank ChatGPT for improving the clarity and conciseness of the sentences in this paper.}}

\usepackage{graphics,color}
\usepackage[outdir=./]{epstopdf}
\usepackage{amsmath}
\usepackage{amssymb}
\usepackage{mathrsfs}
\usepackage{subfigure}
\usepackage[font=small,skip=2pt]{caption}
\usepackage{url}
\usepackage{float}
\usepackage{bbm}
\usepackage{booktabs}
\usepackage{array}
\usepackage[table]{xcolor}
\usepackage{optidef} 

\newtheorem{theorem}{\bf \emph{Theorem}}[section]
\newtheorem{lemma}[theorem]{\bf \emph{Lemma}}

\newtheorem{proposition}[theorem]{Proposition}

\newtheorem{remark}{Remark}
\newtheorem{assumption}[theorem]{Assumption}
\newtheorem{example}{Example}

\newcommand{\Rank}{\operatorname{Rank}}

\newcommand{\transpose}{\mathsf{T}} 

\newcommand*{\QEDB}{\hfill\ensuremath{\blacksquare}}
\newcommand*{\QEDW}{\hfill\ensuremath{\square}}

\newcommand{\trace}{\text{tr}}

\newcommand{\what}{\widehat}
\newcommand{\wtilde}{\widetilde}

\DeclareMathOperator*{\argmin}{arg\,min}

\newcommand{\wtSig}{\widetilde{\Sigma}}

\newcommand{\g}{\beta}

\newcommand{\pl}{{(l)}}

\newcommand{\pzeta}{{(\zeta)}}
\newcommand{\mco}{{\mathcal{O}}}
\newcommand{\mcz}{{\mathcal{Z}}}
\newcommand{\pO}{{(\mco\mco)}}

\newcommand{\pOZ}{{(\mco\mcz)}}

\newcommand{\mbf}{\mathbf}
\newcommand{\bfx}{\mathbf{x}}
\newcommand{\bfy}{\mathbf{y}}
\newcommand{\bfz}{\mathbf{z}}
\newcommand{\bfu}{\mathbf{u}} 
\newcommand{\bs}{\boldsymbol}
\newcommand{\bfA}{\mathbf{A}}
\newcommand{\bfb}{\mathbf{b}}

\newcommand{\bfC}{\mathbf{C}}

\newcommand\oprocendsymbol{\hbox{$\square$}}
\newcommand\oprocend{\relax\ifmmode\else\unskip\hfill\fi\oprocendsymbol}

\renewcommand{\theenumi}{(\roman{enumi})}

\usepackage[linesnumbered,ruled,vlined]{algorithm2e}

\usepackage{tabstackengine}
\setstackEOL{\cr}

\begin{document}
	\maketitle
	
	\thispagestyle{empty} \pagestyle{empty}
	
\begin{abstract} 
In recent papers \cite{al2022observability} and \cite{al2021using} , the second author developed full-state feedback controllers for networked systems to block the observability and controllability of certain remote nodes. In this paper, we build on these control schemes to an interconnected power system with the aims of blocking (i) mode  participation factors and (ii) inter-area mode observability in tie-line power flow measurements. Since participation factors depend on both controllable and observable eigenvectors, the control techniques from the cited works must be carefully tailored to this setting. Our research is motivated by cyber-security concerns in power systems, where an adversary aims to deceive the operator by tampering the system's modal content. We present extensive numerical results on a 3-machine, 9-bus system and a 16-machine, 68-bus system.
\end{abstract}
	
\section{Introduction}
Bulk power systems exhibit a range of dynamic oscillatory behaviors governed by their underlying modes of oscillation. Low frequency oscillations associated with the small-signal (or linear) stability of large power systems are undesirable as they restrict maximum power transfer \cite{rogers2000nature, pal2012coordinated, paserba2001small}. These oscillations are often classified by modes, which represent specific patterns of oscillatory behavior (see Section \ref{sec: prelims} for precise definition). For instance, \emph{local-mode oscillations} arise when generating units within an area oscillate relative to each other. Instead, \emph{inter-area-mode oscillations}, the most crucial and detrimental, occur when large groups of generators in different regions swing against each other.

When disturbances trigger these oscillations in a power system, local or wide-area measurements can capture their effects. It is very essential for system operators to accurately determine which machines or local control units—such as automatic voltage regulators (AVRs) or power system stabilizers (PSS)—are most impacted by the mode. This means understanding whether a machine is being heavily disturbed by the mode or is contributing significantly to it. Such insight allows operators to take prompt corrective actions, like taking a machine or unit offline, to restore stability \cite{machowski2020power}.

A widely used diagnostic tool for identifying which states, such as rotor angle and speed, in a generator (and thus its location) are impacted by an mode is the \emph{mode participation factor}. These factors are defined by the product of the left and right eigenvectors of the system's state matrix corresponding to the mode's eigenvalue. They quantify a state's contribution to a mode. High participation in a controllable state suggests that control actions, such as AVR/PSS tuning, can mitigate the mode, while high participation value in an observable state indicates that the mode can be effectively detected using wide-area measurements (see \cite{sancha1988selective, ghandhari2011stability, abed2000participation, netto2018data})

Recent cyber-physical incidents, including both adversarial and natural disturbances, have demonstrated that power grids are vulnerable to targeted disruptions \cite{dagoumas2019assessing, rajkumar2020cyber, shahkar2020resilient}. Malicious attackers can exploit vulnerabilities in measurement systems and control loops amplifying destabilizing oscillations or misdirecting corrective actions. Similarly, extreme weather events and grid contingencies can induce shifts in inter-area dynamics. These challenges necessitate power engineers to understand the reliability or fragility of monitoring metrics such as participation factors and observing inter-area modes in tie-line power flows. 

We take a first step toward in this direction to understand the extent to which adversaries can change participation factors of any given mode and block observability of inter-area-modes in tie-lie power flow measurements without destabilizing the system. Our approach leverages recently developed state-feedback controllers for blocking observability and controllability in multi-agent networks systems \cite{ al2022observability, al2021using}. 
For an interconnected power system modeled as a linear-time invariant small signal model, we develop a static state-feedback matrix, which, when applied to the power system, results in a closed-loop system with the following attributes: 

\smallskip 
\begin{enumerate}
\item  the modes remain invariant but the participation factors of selective modes are modified, or 
    \item the modes remain invariant but the observable inter-area modes in open-loop system are unobservable. 
\end{enumerate}

We validate the performance our closed-loop controller on a 3-machine, 9-bus system and a 16-machine, 68-bus system. Our results uncover some interesting behavior in altering mode participation factors. In particular, we can block the participation of any desired state in any desired mode, but this comes with a price. The state can now highly participate in unexpected modes. We also discuss the reasons behind this behavior (see Section \ref{sec: simulations}). 

\smallskip

The rest of the paper is as follows. In Section \ref{sec: prelims}, we give a brief introduction to modal analysis in power systems. In Section \ref{sec: problem setup1}, we describe our research problem and in Section \ref{sec: eigenstructure} we present our algorithms for the design of blocking controllers. In Section \ref{sec: simulations}, we provide simulation results and conclude the paper in Section \ref{sec: conclusions}. 


\section{A Primer on Power System Dynamics}\label{sec: prelims}
Since our interest is in the time scale where low-frequency oscillations manifest in an interconnected power system, we consider a small signal stable model that can be obtained by linearizing non-linear differential algebraic equations \cite{kundur2007power}:
\begin{subequations}\label{eq: DAE power equations}
\begin{align}
    \dot{x}(t) &= f(x(t),z(t),u(t)) \label{eq: DAE1}\\
    0 &= g(x(t),z(t),u(t)) \label{eq: DAE2}\\
    y(t) &= h(x(t),z(t),u(t)) \label{eq: DAE3}
\end{align}
\end{subequations}
where $x$, $z$, and $u$ are the state, the algebraic, and the input vectors, resp.; and $f$, $g$, and $h$ are vector-valued functions. 

The state $x(t)\in \mathbb{R}^n$ comprises rotor angles, rotor speeds, flux linkages, and the states of local dynamic components such as exciters, automatic voltage regulators (AVRs), power system stabilizers (PSS), and governors for all machines. The vector $z(t)\in \mathbb{R}^r$ comprises algebraic variables such as bus phase angles and voltage magnitudes, and $d$- and $q$-axes currents in the stator equations. The control input $u(t)\in \mathbb{R}^q$ typically include mechanical power, set points of AVRs, etc.

Linearizing around the operating point $(x_0,u_0,y_0)$ and eliminating the algebraic variable $z(t)$ from \eqref{eq: DAE power equations} yields us
\begin{equation}\label{eq: linear systems}
\left\{\begin{array}{l}
\Delta \dot{{x}}(t)=\mathbf{A} \Delta {x}(t)+\mathbf{B}{u}(t) \\
\Delta y(t)=\mathbf{C} \Delta {x}(t),
\end{array}\right.
\end{equation}
where $\Delta$ signifies that the dynamics are deviations from the operating point, while the state-space matrices $(\mathbf{A}, \mathbf{B}, \mathbf{C})$ are functions of the operating point. Moreover, these matrices can be obtained either analytically or numerically. While the numerical approach is preferred for large-scale simulations, the analytical method provides insight into the specialized structure of the state-space matrices; see Section \ref{sec: simulations}. 

\subsection{Modal Analysis}
We give a brief overview of modal analysis. For a detailed exposition, see \cite{machowski2020power}. The eigenvalue of $\mathbf{A}$ is called the \emph{mode}. Let $\Lambda$ be the diagonal matrix of modes and they are distinct. Then $\mathbf{A}V = V \Lambda$ and $W^\transpose \mathbf{A}=\Lambda W^\transpose$, where the columns of the matrix $V$ are the right eigenvectors and the rows of the matrix $W^\transpose$ are the left eigenvectors. 
These matrices, also known as modal matrices, satisfy the relation $W^\transpose V=I$ and diagonalize the matrix $\mathbf{A}$; that is, $\Lambda=W^\transpose\mathbf{A}V$. 

Define the \emph{modal transform} $z=W^\transpose\Delta x$. Then $\dot{z}=\Lambda z$ is the modal form for the dynamics given in \eqref{eq: linear systems} with no input. Further, the \emph{modal variable} $z_i(t)$ can be expressed as 
\begin{align*}
    z_i(t)=w_{i1}\Delta x_{1}(t)+w_{i2}\Delta x_{2}(t)+\ldots+w_{in}\Delta x_{n}(t), 
\end{align*}
where $w_{ij}$ is the $j$-th entry of the $i$-th row vector of $W$ (that is, the left eigenvector of the mode $i$). And $w_{ij}$ determines the \emph{magnitude and phase of the share of $x_j(t)$ in the activity of the mode $z_i(t)$}. Hence, a large $w_{ij}$ means controlling state $j$ will strongly excite mode $i$. 

Similarly, the state variable $x_k(t)$ can be written as
\begin{align*}
    \Delta x_k(t)=v_{k1}z_{1}(t)+v_{k2}z_{2}(t)+\ldots+v_{kn}z_{n}(t), 
\end{align*}
where $v_{ki}$ is the $k$-th entry of the $i$-th column of $W$ (the right eigenvector of the mode $i$), and it determines the \emph{magnitude and phase of the share of $z_i(t)$ in the activity of $x_k(t)$}. Hence, a large $v_{ki}$ means that state $k$ is highly active in mode $i$. 



\subsection{Participation Factors}
Participation factors quantify the contribution of each state variable to a specific mode of the system. They are used as a diagnostic tool to decide which machines are contributing to a particular oscillation and also for placing a damping controller. Mathematically, participation factor for state $x_k(t)$ in modal variable $z_i(t)$ is given by the product of the $k$-th element of the $i$-th left and right eigenvectors: $p_{ki}=w_{ik}v_{ki}$. All factors can be computed at once using the participation matrix \( P = V \circ W^\transpose \), where $\circ$ is the Hadamard product. 

The interpretation of these factors is clear if we consider the impact of the \( i \)-th swing mode (an oscillatory mode that arises due to the interaction of rotor speeds and angles) on the \( k \)-th state associated with machine \( M1 \). A larger \( p_{ki} \) means that \( M1 \) either significantly swings in mode \( i \) or strongly excites it. Thus, \( M1 \) has high participation in mode \( i \).  

\subsection{Inter-Area Mode Observability in Tie-Line Power Flows}
In modal coordinates the entries of the measurement vector $\Delta y(t)$ can be written as 
\begin{align}
    y_j(t)=\sum_{i=1}^no_{ji}z_i(t) \quad j=1,\ldots,p, 
\end{align}
where $o_{ji}=c^\transpose_jv_i$ and $c_j^\transpose$ is the $j$-th row of $\mathbf{C}$. Thus, if $o_{ji}=0$ for all $j$ the $i$-th mode is unobservable in measurements. In our simulation study we focus on a set of $y_j$ measurments rather than one single measurement. 

As discussed in the introduction, two types of modes are crucial in power system. Local modes (0.7–2 Hz) involve a single generator or a group oscillating against the system, while inter-area modes (0.1–1 Hz) occur when groups of machines in different regions swing out of phase, imposing significant limits on tie-line power flows. Tie-lines connect different regions of the grid. Inter-area modes are unobserved in local measurements of generators and are clearly seen in global measurements like tie-line power flows.

\section{Research question}\label{sec: problem setup1}
We make two assumptions. First, a centralized authority can jointly actuate all available control inputs (e.g., AVR set points or PSS inputs). This assumption is well-justified and commonly used to mitigate inter-area oscillations through coordination of multiple local controllers \cite{abdalla1984coordinated}. Second, all state variables are measurable. While  restrictive, this assumption is reasonable when using a reduced-order or aggregated model \cite{chow2013power}. With these assumptions in place, we address: 
\begin{enumerate}
    \item for a given set of modes (local or inter-area), design a full-state feedback controller to block participation of a few select states in those modes.
    \item design a full-state feedback controller to block observability of inter-area mode in tie-line flows. 
\end{enumerate}

Mathematically, the first question is to enforce $p_{ki}=0$ for a specific state $k$ and a mode $i$, effectively eliminating that state's participation in the mode. The second objective aims to set $o_{ji}=0$ for all inter-area modes $i$ in the $j$-th tie-line power flow measurements.


\section{Surgical Eigenstructure Assignment}\label{sec: eigenstructure}

Using our previous results on the surgical eigenstructure assignment problem \cite{al2022observability,al2021using,al2019observability} we present algorithms to compute a full-state feedback control law $u(t)=\mathbf{F}x(t)$ to (i) shape participation factors and (ii) enforce unobservability in tie-line power flow measurements for given modes.

\subsection{Blocking Participation Factors}
Suppose that we want $m$ states to not participate in the mode $\lambda_i$ of the state matrix $\mathbf{A}$. Then our algorithm outputs a state feedback matrix $\mathbf{F}$ such that $p_{k_1i}=p_{k_2i}=\cdots =p_{k_mi}=0$
for the mode $\lambda_i$ associated with the closed-loop matrix $\mathbf{(A+BF)}$. We emphasize that our state feedback keeps the mode $\lambda_i$ invariant while modifying the $k_1$-th, $k_2$-th, $\cdots$, $k_m$-th entries of the open-loop right eigenvector $\mathbf{v}_i$ to become zero. These zero entries in the eigenvector will enforce zero participation of corresponding states in the mode $\lambda_i$. Below we present our algorithm for blocking participation factors in a complex conjugate mode pair. For more details, see \cite{al2022observability}. 

\medskip 


\noindent \textbf{Algorithm 1:} 

\smallskip 
1) Select a pair of complex conjugate mode $\lambda_i$ and ${\lambda}_{i+1}$ of $\mathbf{A}$ and its associated eigenvector $\mathbf{v}_i$ and $\mathbf{v}_{i+1}$ where $i \in \{1, \hdots, n-1\}$. Note, ${\lambda}_{i+1}=\bar{\lambda}_i$ and $\mathbf{v}_{i+1}=\bar{\mathbf{v}}_i$.

\smallskip 
2) Compute $\mathbf{N}(\lambda_i) \in \mathbb{C}^{(n+q) \times q}$, whose columns are linearly independent and span the null space of the matrix $\mathbf{S}(\lambda_i)=[(\mathbf{A}-\lambda_i\mathbf{I}_n)~~ \mathbf{B}]$.  Consider the partition 
$\mathbf{N}(\lambda_i)=[\mathbf{N}_1(\lambda_i)^\transpose ~ \mathbf{N}_2(\lambda_i)^\transpose]^\transpose$, where  $\mathbf{N}_1(\lambda_i) \in \mathbb{C}^{n \times q}$ and  $\mathbf{N}_2(\lambda_i) \in \mathbb{C}^{q \times q}$. Therefore $\mathbf{N}_1(\lambda_i)$ and $\mathbf{N}_2(\lambda_i)$ satisfy:
\begin{eqnarray}
[(\mathbf{A}-\lambda_i~\mathbf{I}_n)~~ \mathbf{B}] \left[\begin{array}{c}
                                         \mathbf{N}_1(\lambda_i)   \\
                                         \mathbf{N}_2(\lambda_i)
                                          \end{array}\right]  = \mathbf{0}.  \label{eq4}
\end{eqnarray}

3) Construct the matrix $\mathbf{N}_3(\lambda_i) \in \mathbb{C}^{m \times q}$ by selecting $k_1$-th, $k_2$-th, ..., $k_m$-th rows of $\mathbf{N}_1(\lambda_i)$. Then find a vector $\mathbf{h}_i \neq \mathbf{0}$ which lies in the null space of $\mathbf{N}_3(\lambda_i)$, i.e. which satisfies:  
\begin{eqnarray}
\mathbf{N}_3(\lambda_i)~ \mathbf{h}_i &=& \mathbf{0}. \label{eq2}
\end{eqnarray}

\smallskip 
4) Compute the vectors $\hat{\mathbf{v}}_i$ and $\mathbf{z}_i$ as:
\begin{eqnarray}
\hat{\mathbf{v}}_i &=& \mathbf{N}_1(\lambda_i)~ \mathbf{h}_i \label{eq1} \\
\mathbf{z}_i &=& \mathbf{N}_2(\lambda_i)~ \mathbf{h}_i. \label{eq5}
\end{eqnarray}
Choose $\hat{\mathbf{v}}_{i+1}=\bar{\hat{\mathbf{v}}}_i$ and $\mathbf{z}_{i+1}=\bar{\mathbf{z}}_i$.

\smallskip 
5) Construct $\hat{\mathbf{V}}$ from $\mathbf{V}$ by replacing the columns containing $\mathbf{v}_i$ and $\mathbf{v}_{i+1}$ with $\hat{\mathbf{v}}_i$ and $\hat{\mathbf{v}}_{i+1}$. In the same manner construct $\mathbf{Z}$ from $\mathbf{Z}_0=\mathbf{0}\in \mathbb{R}^{q\times n}$ by replacing the corresponding columns of $\mathbf{Z}_0$ with $\mathbf{z}_i$ and  $\mathbf{z}_{i+1}$ obtained in Step 4.  Therefore, $\hat{\mathbf{V}}=[ \mathbf{v}_1 \cdots \mathbf{v}_{i-1}~ \hat{\mathbf{v}}_i~\hat{\mathbf{v}}_{i+1}~\mathbf{v}_{i+2}~ \cdots \mathbf{v}_n]$ and $\mathbf{Z}= [\mathbf{0} \cdots \mathbf{0} ~ \mathbf{z}_i~ \mathbf{z}_{i+1}~ \mathbf{0} \cdots \mathbf{0}]$.

6) Finally the gain matrix $\mathbf{F}$ is:
\begin{eqnarray}
 \mathbf{F} = \mathbf{Z}~\hat{\mathbf{V}}^{-1}.  \label{eq3a}
\end{eqnarray}

We assume that the matrix $\hat{\mathbf{V}}$ constructed in Step (5) is invertible. If not, additional steps must be taken to guarantee invertibility. We refer the reader to \cite{al2022observability} for details on this. If $m+2\leq q$, then according to \cite{al2022observability}, Algorithm 1 is guaranteed to find the desired $\mathbf{F}$ to block participation.  

\subsection{Blocking Observability for Tie-Line Power Flow}
Blocking observability for the tie-line power flow $y(t)$ corresponds to modifying open-loop eigenvector $\mathbf{v}_i$ of $\mathbf{A}$ to $\hat{\mathbf{v}}_i$ so that $\mathbf{C} \mathbf{v}_i=\mathbf{0}$ where $\hat{\mathbf{v}}_i$ is the modified eigenvector of $\mathbf{(A+BF)}$ corresponding to the mode $\lambda_i$. This operation will render the mode $\lambda_i$ to become unobservable according to the PBH test \cite{rugh1996linear}. We present the algorithm for blocking observability in a pair of complex conjugate modes.

\medskip 
\noindent \textbf{Algorithm 2:} 

1) Same as Step 1 in Algorithm 1.

\smallskip 
2) Same as Step 2 in Algorithm 1.

\smallskip 
3) Find the matrix $\mathbf{M}_1$ whose columns form a basis of the null space of $\mathbf{C}$. Then construct the matrix $\mathbf{M}_2(\lambda_i)= [\mathbf{N}_1(\lambda_i)~~\mathbf{M}_1]$. Find the matrix $\mathbf{M}_3(\lambda_i)$ whose columns form a basis of the null space of $\mathbf{M}_2(\lambda_i)$. Then $\mathbf{h}_i$ is obtained by the first $q$ entries of any columns of $\mathbf{M}_3(\lambda_i$).

\smallskip 
4-6) Steps 4-6 are the same as for Algorithm 1.

As before, we assume that the matrix $\hat{\mathbf{V}}$ constructed in Step (5) is invertible. If not, additional steps must be taken to guarantee invertibility (see \cite{al2022observability}). Using the arguments in \cite{al2022observability}, we can show that if $rank(\mathbf{C})+2\leq q$, then Algorithm 2 is guaranteed to find the desired $\mathbf{F}$ to block observability. 

\medskip 
An important point here is that Step 3 of Algorithm 2 is different than the algorithms in \cite{al2022observability} as here we allow a more general structure for $\mathbf{C}$ matrix. In our algorithm, we compute $\hat{\mathbf{v}}_i$ in a way so that it lies in the intersection of the null space of $\mathbf{C}$ and assignable eigenvector space $\mathbf{N}_1(\lambda_i)$. We further note that our algorithms maintain the eigenvalues and right eigenvectors except $\hat{\mathbf{v}}_i$  and $\hat{\mathbf{v}}_{i+1}$, therefore they can be used repeatedly to block observability or participation in multiple modes.

\section{Participation and inter-area-observability blocking in power systems}\label{sec: simulations}  
We use Algorithms 1 and 2 in Section \ref{sec: eigenstructure} to block selected mode participation factors and specific inter-area modes in the tie-line power flows of a (i) 3-machine, 9-bus system and (ii) a 16-machine, 68-bus system. For both systems, at an appropriate operating point, eliminating algebraic variables yields the state-space models given in \eqref{eq: linear systems}. To highlight the impact of our control method, for the three machine system, we work with the analytical Heffron-Phillips model.

\subsection{Case 1: 3 machine, 9 bus power system:} We adopt the presentation in \cite{wang2016analysis}, assuming a fourth-order model for each machine equipped with a local PSS unit. Let $\Delta\omega=[\Delta\omega_1~ \Delta\omega_2~\Delta\omega_3]^\transpose$ be the rotor speed in rad/sec. Similarly, define the state vectors: rotor angle $\Delta \delta$ (rad), $q$-axis transient EMF $\Delta E'_q$, and field excitation voltage $E'_{fd}$. 

Define $\Delta x=[\Delta\delta^\transpose~\Delta\omega^\transpose~ (\Delta E'_q)^\transpose~ (\Delta E'_{fd})^\transpose]^\transpose \in \mathbb{R}^{12}$ and $\Delta u_\text{PSS}=[\Delta u_\text{PSS}~ \Delta u_\text{PSS}~ \Delta u_\text{PSS}]^\transpose \in \mathbb{R}^{3}$ be the output of PSS. Then, the state-space matrices take the form: 
\begin{align*}
\mathbf{A}\!=\!\left[\begin{array}{cccc}
\mathbf{0} & {\omega}_0 \mathbf{I} & \mathbf{0} & \mathbf{0} \\
-\mathbf{M}^{-1} \mathbf{K}_1 & -\mathbf{M}^{-1} \mathbf{D} & -\mathbf{M}^{-1} \mathbf{K}_{\mathbf{2}} & \mathbf{0} \\
-\mathbf{T}_{\mathbf{d} \mathbf{0}}^{-1} \mathbf{K}_{\mathbf{4}} & \mathbf{0}  & -\mathbf{T}_{d 0}^{-1} \mathbf{K}_{\mathbf{3}} & \mathbf{T}_{\mathrm{d} \mathbf{0} }^{-1} \\
-\mathbf{T}_{\mathbf{A}}^{-1} \mathbf{K}_5 \mathbf{K}_{\mathbf{A}} & \mathbf{0} & -\mathbf{T}_{\mathbf{A}}^{-1} \mathbf{K}_5 \mathbf{K}_{\mathbf{A}} & -\mathbf{T}_{\mathbf{A}}^{-1}
\end{array}\right]
\end{align*}
\begin{align*}
\mathbf{B}=\left[\begin{array}{c}
\mathbf{0} \\
\mathbf{0} \\
\mathbf{0} \\
\mathbf{T}_{\mathbf{A}}^{-1} \mathbf{K}_{\mathbf{A}}
\end{array}\right], 
\end{align*}
where $\mathbf{0}$ and $\mathbf{I}$ are zero and identify matrices of appropriate dimensions, and $\omega_0$ is the reference speed. The structure and numerical values of block matrices in $\mathbf{A}$ and $\mathbf{B}$ can be found in \cite[Pages 196 and 220]{wang2016analysis}.



Let $P=[P_1~P_2~ P_3]^\transpose \in \mathbb{R}^{3}$ be the vector of net electrical power injections at all machines. The deviations in net power injection then obey $\Delta P=\mathbf{K}_1\Delta \delta + \mathbf{K}_2 \Delta E'_q$. Then the tie-line power flows between all three machines are given by: 
\begin{align*}
P_\text{Tie-line}\triangleq \begin{bmatrix}
     P_{12}\\
     P_{13}\\
     P_{23}
 \end{bmatrix}=\underbrace{\begin{bmatrix}
     1 & -1 & 0\\
     1 & 0 & -1\\
     0 & 1 & -1
 \end{bmatrix}}_{\triangleq \widetilde{\boldsymbol{C}}} \begin{bmatrix}
     \Delta P_{1}\\
     \Delta P_{2}\\
     \Delta P_{3}
 \end{bmatrix},   
\end{align*}
which can be expressed in terms of the state vector $\Delta x$ as 
\begin{align*}
    P_\text{Tie-line}=\underbrace{\widetilde{\boldsymbol{C}}\begin{bmatrix}
         \mbf{K}_1 & \mathbf{0} & \mbf{K}_2 & \mathbf{0}
    \end{bmatrix}}_{\triangleq \mathbf{C}}\Delta x
\end{align*}

\begin{figure*}[t] 
    \centering
    \includegraphics[width=0.95\linewidth]{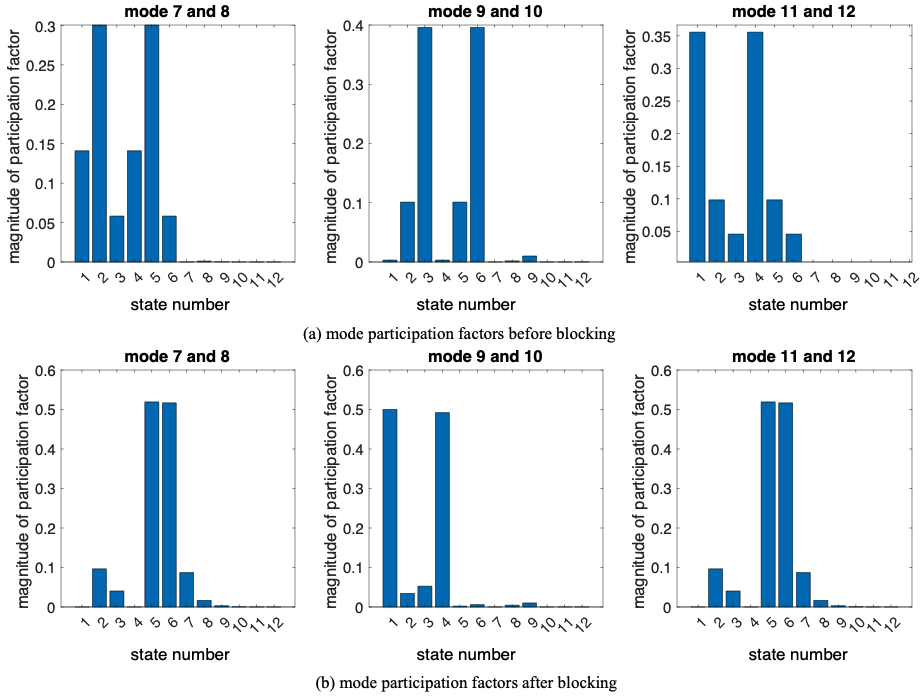} 
    \caption{Blocking participation of states 1 and 4 of of 3 machine 9 bus power system using Algorithm 1. Refer main text for details.}
    \label{fig: 3machine_PFs}
\end{figure*}

We now use Algorithm 1 to block the participation of state 1 (the rotor angle $\Delta \delta_1$ of machine 1) and state 4 (speed $\Delta \omega_1$ of machine 1) in the highly participating modes. In Fig.~\ref{fig: 3machine_PFs}(a), we see that these states highly participate in mode 7, mode 8 (the conjugate of mode 7), mode 11, and finally, mode 12 (the conjugate of mode 11). Fig.~\ref{fig: 3machine_PFs}(b) shows the magnitudes of mode participation factors after employing our feedback controller. Due to space constraints, we do not display participation factors for other modes, as states 1 and 4 do not highly participate in them.

\smallskip 
In Fig.~\ref{fig: 3machine_PFs}(b), we see that states 1 and 4 no longer participate in modes 7, 8, 11, and 12. However, interestingly, these states participate in modes 9 and 10. This behavior arises because our state feedback modifies only the right eigenvectors of modes 7, 8, 11, and 12 while leaving the remaining eigenvectors unchanged, whereas the left eigenstructure undergoes a complete transformation. Since participation factors are computed using the Hadamard product of right and left eigenvectors, the new participation factors may either increase or decrease. In future, we will develop feedback controllers that also constrain changes to the left eigenvectors, providing more precise control over participation factors. We finally note that the participation factors of all other modes (i.e., modes 1-6) are almost the same as before.

\begin{figure*}[htpb]
    \centering
    \includegraphics[width=0.8\linewidth]{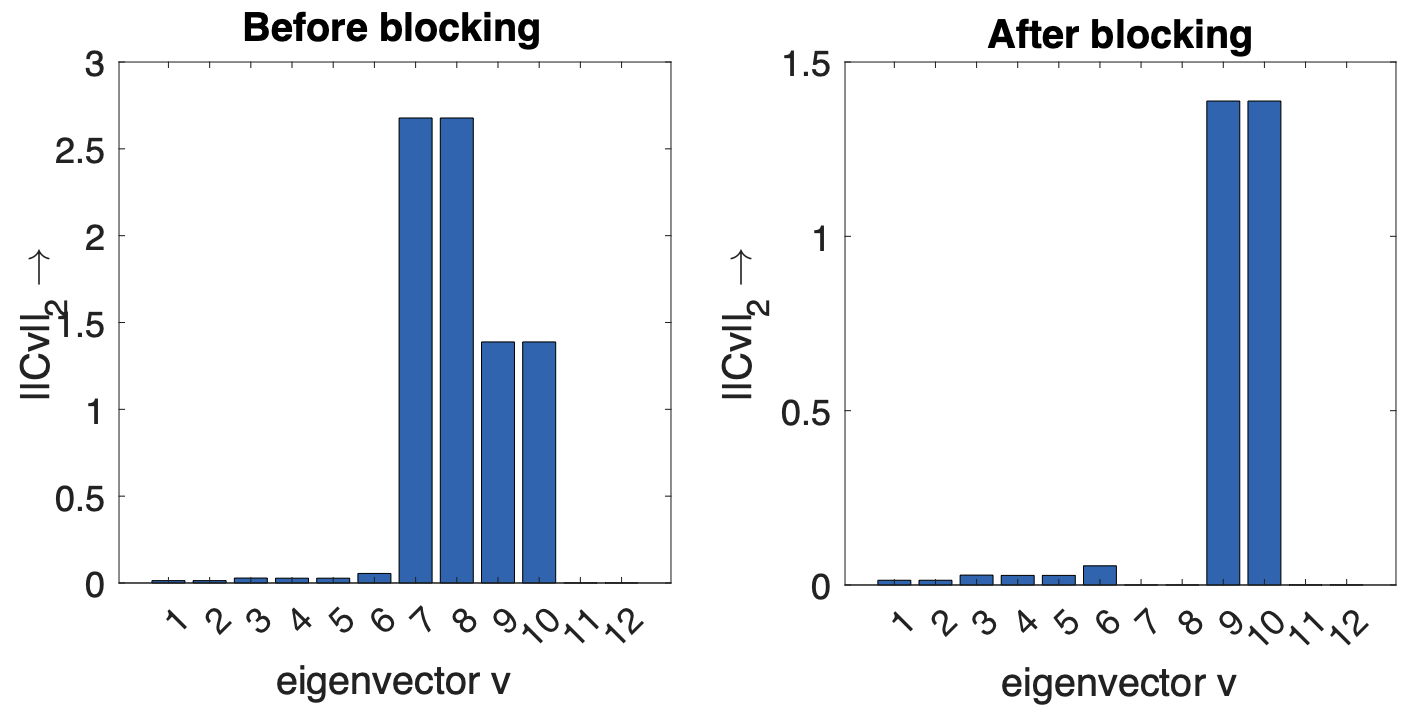}
    \caption{Norm of $\mathbf{C}\mathbf{v}_i$ for each eigenvector $\mathbf{v}_i$ of 3 machine 9 bus power system before (left) and after (right) inter-area blocking.}
    \label{fig:M3}   
\end{figure*}

\smallskip 
We block the observability for the measurement of tie-line power flow using Algorithm 2. Specifically, we enforce unobservability on mode 7 and mode 8 (conjugate of mode 7) for output matrix $\mathbf{C}$ associated with tie-line power flow measurements. Recall that according to the PBH condition, the mode $\lambda_i$ is unobservable if $\mathbf{C} \mathbf{v}_i=\mathbf{0}$ where $\mathbf{v}_i$ is the eigenvector corresponding to mode $\lambda_i$. To depict unobservability, we thereby present the norm of $\mathbf{C}\mathbf{v}_i$ for all the eigenvectors in Fig. \ref{fig:M3}. As shown in the figure, the norms of the eigenvectors corresponding to modes 7 and 8 are zero, indicating that unobservability is successful on these modes.

\smallskip 
\subsection{Case 2: 16 machine, 68 bus power system:} 
Each machine has ten states, and hence the total number of states is $10 \times 16 = 160$. The states are grouped such that those associated with each generator are ordered together: 
\begin{align*}
\mathbf{x}^T\!=\!(x^{(1)}_{1} \dots x^{(1)}_{10} \mid x^{(2)}_{1} \dots x^{(2)}_{10} \mid \dots \mid x^{(16)}_{1} \dots x^{(16)}_{10}). 
\end{align*}
Thus, for example, $x^{(1)}_{5}$ represents the fifth state of generator 1. For the \( i \)-th generator, the state vector components are:  
\( x^{(i)}_{1} = \delta \) (rotor angle),  
\( x^{(i)}_{2} = \omega \) (speed),  
\( x^{(i)}_{3} = E'_q \) (q-axis transient internal voltage),  
\( x^{(i)}_{4} = \psi^{"}_d \) (d-axis flux),  
\( x^{(i)}_{5} = E'_d \) (d-axis transient internal voltage),  
\( x^{(i)}_{6} = \psi^{"}_q \) (q-axis flux),  
\( x^{(i)}_{7} = E_{fd} \) (generator field voltage),  
\( x^{(i)}_{8}, x^{(i)}_{9}, x^{(i)}_{10} \) (internal states of the PSS). We use the power system toolbox \cite{sauer2017power} to obtain state-space matrices with dimensions: $\mbf{A} \in \mathbb{R}^{160\times 160}$ and $\mbf{B}\in \mathbb{R}^{160\times 16}$. The inputs are exciters' field voltages.

\begin{figure*}[htpb]
    \centering
    \includegraphics[width=0.8\linewidth]{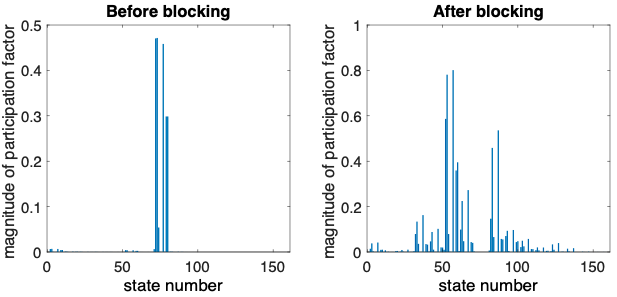}
    \caption{Participation factors of the modes 6 (and associate conjugate mode 7) for each state of 16 machine 68 bus power system.}
    \label{fig:M4}   
\end{figure*}

\smallskip 
We block the participation of all ten states of a machine in a particular mode.  We chose mode 6 and mode 7 (conjugate of mode 6). This conjugate mode pair is a local mode because states of machine M7 (for one line diagram of 16 bus system, see \cite{singh2013report}) participate in them. Now, we block the participation of the ten states of M7 to this conjugate mode pair using Algorithm 1. As displayed on the right-hand side of Fig.~\ref{fig:M4}, states associated with machine 7 do not participate anymore in modes 6 and 7. However, several states of other machines now highly participate in this mode, implying non-locality of the mode. Finally, to avoid repetition, we do not present results on blocking inter-area mode observability.  

\section{Concluding Remarks}\label{sec: conclusions}
This paper proposes full-state feedback control schemes to modify and suppress mode participation factors and inter-area mode observability in interconnected power systems. The control algorithm leverages the classical eigenstructure assignment method. Through detailed simulations on a 3-machine, 9-bus system and a 16-machine, 68-bus system, we demonstrate that our approach effectively reveals fundamental limits in power systems relying on modal analysis.

However, a drawback of this study is that full-state feedback is impractical in modern power systems. This highlights the need for future research on developing output feedback or decentralized control strategies for blocking mode participation factors and observability.

\bibliographystyle{unsrt}
\bibliography{BIB.bib}

\begin{thebibliography}{10}

\bibitem{al2022observability}
Abdullah Al~Maruf and Sandip Roy.
\newblock Observability-blocking control using sparser and regional feedback for network synchronization processes.
\newblock {\em Automatica}, 146:110586, 2022.

\bibitem{al2021using}
Abdullah Al~Maruf and Sandip Roy.
\newblock Using feedback to block controllability at remote nodes in network synchronization processes.
\newblock In {\em 2021 American Control Conference (ACC)}, pages 2473--2478. IEEE, 2021.

\bibitem{rogers2000nature}
Graham Rogers.
\newblock The nature of power system oscillations.
\newblock In {\em Power System Oscillations}, pages 7--30. Springer, 2000.

\bibitem{pal2012coordinated}
Anamitra Pal.
\newblock {\em Coordinated Control of Inter-area Oscillations using SMA and LMI}.
\newblock PhD thesis, Virginia Tech, 2012.

\bibitem{paserba2001small}
John Paserba, Juan Sanchez-Gasca, Prabha Kundur, Einar Larsen, and Charles Concordia.
\newblock Small signal stability and power system oscillations.
\newblock {\em Electric power engineering handbook}, pages 20--34, 2001.

\bibitem{machowski2020power}
Jan Machowski, Zbigniew Lubosny, Janusz~W Bialek, and James~R Bumby.
\newblock {\em Power system dynamics: stability and control}.
\newblock John Wiley \& Sons, 2020.

\bibitem{sancha1988selective}
JL~Sancha and IJ~Perez-Arriaga.
\newblock Selective modal analysis of power system oscillatory instability.
\newblock {\em IEEE Transactions on Power Systems}, 3(2):429--438, 1988.

\bibitem{ghandhari2011stability}
Mehrdad Ghandhari.
\newblock Stability of power systems.
\newblock {\em Electric Power Systems, Royal Institute of Technology, Stockholm, Sweden}, 2011.

\bibitem{abed2000participation}
Eyad~H Abed, David Lindsay, and Wael~A Hashlamoun.
\newblock On participation factors for linear systems.
\newblock {\em Automatica}, 36(10):1489--1496, 2000.

\bibitem{netto2018data}
Marcos Netto, Yoshihiko Susuki, and Lamine Mili.
\newblock Data-driven participation factors for nonlinear systems based on koopman mode decomposition.
\newblock {\em IEEE control systems letters}, 3(1):198--203, 2018.

\bibitem{dagoumas2019assessing}
Athanasios Dagoumas.
\newblock Assessing the impact of cybersecurity attacks on power systems.
\newblock {\em Energies}, 12(4):725, 2019.

\bibitem{rajkumar2020cyber}
Vetrivel~Subramaniam Rajkumar, Marko Tealane, Alexandru {\c{S}}tefanov, Alfan Presekal, and Peter Palensky.
\newblock Cyber attacks on power system automation and protection and impact analysis.
\newblock In {\em 2020 IEEE PES Innovative Smart Grid Technologies Europe (ISGT-Europe)}, pages 247--254. IEEE, 2020.

\bibitem{shahkar2020resilient}
Shahram Shahkar and Khashayar Khorasani.
\newblock A resilient control against time-delay switch and denial of service cyber attacks on load frequency control of distributed power systems.
\newblock In {\em 2020 IEEE Conference on Control Technology and Applications (CCTA)}, pages 718--725. IEEE, 2020.

\bibitem{kundur2007power}
Prabha Kundur.
\newblock Power system stability.
\newblock {\em Power system stability and control}, 10:7--1, 2007.

\bibitem{abdalla1984coordinated}
Omar~H Abdalla, SA~Hassan, and NT~Tweig.
\newblock Coordinated stabilization of a multimachine power system.
\newblock {\em IEEE Transactions on power Apparatus and Systems}, (3):483--494, 1984.

\bibitem{chow2013power}
Joe~H Chow.
\newblock {\em Power system coherency and model reduction}, volume~84.
\newblock Springer, 2013.

\bibitem{al2019observability}
Abdullah Al~Maruf and Sandip Roy.
\newblock Observability-blocking controllers for network synchronization processes.
\newblock In {\em 2019 American Control Conference (ACC)}, pages 2066--2071. IEEE, 2019.

\bibitem{rugh1996linear}
Wilson~J Rugh.
\newblock {\em Linear system theory}.
\newblock Prentice-Hall, Inc., 1996.

\bibitem{wang2016analysis}
Haifeng Wang, Wenjuan Du, et~al.
\newblock {\em Analysis and damping control of power system low-frequency oscillations}, volume~1.
\newblock Springer, 2016.

\bibitem{sauer2017power}
Peter~W Sauer, MA~Pai, and Joe~H Chow.
\newblock Power system toolbox.
\newblock 2017.

\bibitem{singh2013report}
Abhinav~Kumar Singh and Bikash~C Pal.
\newblock Report on the 68-bus, 16-machine, 5-area system.
\newblock 2013.

\end{thebibliography}







\end{document}